\documentclass{elsart3-1}

\usepackage{amssymb}
\usepackage{amsmath,latexsym,graphicx}
\usepackage[english,francais]{babel}

\newtheorem{theorem}{Theorem}[section]
\newtheorem{lemma}[theorem]{Lemma}
\newtheorem{e-proposition}[theorem]{Proposition}
\newtheorem{corollary}[theorem]{Corollary}
\newtheorem{e-definition}[theorem]{Definition\rm}
\newtheorem{remark}{\it Remark\/}

\newtheorem{theoreme}{Th\`eor\`eme}[section]

\newtheorem{proposition}[theoreme]{Proposition}

\renewcommand{\theequation}{\arabic{equation}}
\DeclareMathOperator{\sech}{sech}

\def\og{\leavevmode\raise.3ex\hbox{$\scriptscriptstyle\langle\!\langle$~}}
\def\fg{\leavevmode\raise.3ex\hbox{~$\!\scriptscriptstyle\,\rangle_0\!\rangle_0$}}

\makeatletter
\newcommand{\defeq}{\mathrel{\rlap{\raisebox{0.3ex}{$\m@th\cdot$}}\raisebox{-0.3ex}{$\m@th\cdot$}}=}
\makeatother

\newcommand{\dk}{\,\mathrm{d}k}
\newcommand{\dx}{\,\mathrm{d}x}

\renewcommand{\e}{\mathrm{e}}
\renewcommand{\i}{\mathrm{i}}

\newcommand{\eqn}[1]{(\ref{#1})}

\newcommand{\nn}{|{\mskip-2mu}|{\mskip-2mu}|}

\journal{the Acad\`emie des sciences}
\begin{document}
\centerline{}

\vspace{-4cm}
\begin{frontmatter}


\selectlanguage{english}
\title{Variational existence theory for hydroelastic solitary waves}


\selectlanguage{english}
\author[sb,leic]{Mark D. Groves},
\ead{groves@math.uni-sb.de}
\author[sb]{Benedikt Hewer},
\author[lund]{Erik Wahl\'{e}n}

\address[sb]{Fachrichtung Mathematik, Universit\"{a}t des Saarlandes, Postfach 151150, 66041 Saarbr\"{u}cken, Germany}
\address[leic]{Department of Mathematical Sciences, Loughborough University, Loughborough, Leics, LE11 3TU, UK}
\address[lund]{Centre for Mathematical Sciences, Lund University, PO Box 118, 22100 Lund, Sweden}


\begin{center}
{\small Received *****; accepted after revision +++++\\
Presented by *****}
\end{center}
\vspace{-0.5cm}\begin{abstract}
\selectlanguage{english}
This paper presents an existence theory for
solitary waves at the interface between a thin ice sheet (modelled using the Cosserat theory of hyperelastic shells)
and an ideal
fluid (of finite depth and in irrotational motion) for sufficiently large values of a dimensionless
parameter $\gamma$. We establish
the existence of a minimiser of the wave energy ${\mathcal E}$
subject to the constraint ${\mathcal I}=2\mu$, where ${\mathcal I}$ is the
horizontal impulse and $0< \mu \ll 1$, and show that the solitary waves detected by
our variational method converge (after an appropriate rescaling)
to solutions of the nonlinear
Schr\"{o}dinger equation with cubic focussing nonlinearity as $\mu \downarrow 0$.
\vskip 0.5\baselineskip

\selectlanguage{francais}
\noindent{\bf R\'{e}sum\'{e}} \vskip 0.5\baselineskip \noindent
{\bf Une th\'{e}orie variationnelle d'existence d'ondes solitaires hydro\'{e}lastiques.}
Cette note pr\'{e}sente une th\'{e}orie d'existence d'ondes solitaires \`{a} l'interface entre une couche de glace mince (mod\'{e}lis\'{e}e par la th\'{e}orie des coques hyper\'{e}lastiques de Cosserat) et un fluide parfait (de profondeur finie et irrotationnel), pour des valeurs suffisamment grandes d'un param\`{e}tre sans dimension $\gamma$. Nous montrons l'existence d'un minimiseur de l'\'{e}nergie ${\mathcal E}$ de l'onde sous la contrainte ${\mathcal I}=2\mu$, o\`{u} ${\mathcal I}$ repr\'{e}sente l'impulsion horizontale et $0< \mu \ll 1$. Nous d\'{e}montrons que les ondes solitaires trouv\'{e}es par notre m\'{e}thode variationnelle convergent (apr\`{e}s un changement d'\'{e}chelle appropri\'{e}) vers des solutions de l'\'{e}quation de Schr\"{o}dinger cubique focalisante, lorsque $\mu \downarrow 0$.

\end{abstract}
\end{frontmatter}

\selectlanguage{english}

\addtolength{\abovedisplayskip}{-2.5pt}
\addtolength{\belowdisplayskip}{-2.5pt}

\vspace{-1.25cm}\section{Introduction}

\vspace{-0.25cm}\subsection{The hydrodynamic problem}

\vspace{-0.2cm}In this article we consider the two-dimensional irrotational flow of a perfect fluid beneath a thin ice sheet
modelled using the Cosserat theory of hyperelastic shells (Plotnikov and Toland \cite{PlotnikovToland11}).
The fluid is bounded below by a rigid horizontal bottom $\{y=0\}$ and above by a free surface
$\{y=h+\eta(x,t)\}$; there is no cavitation between this surface and the ice sheet.
The mathematical problem is to find an Eulerian velocity potential $\phi$ which satisfies the equations
\begin{align}
& \parbox{6.5cm}{$\phi_{xx}+\phi_{yy}=0$,} 0<y<1+\eta, \label{HW1} \\
& \parbox{6.5cm}{$\phi_y = 0$,} y=0, \label{HW2} \\
& \parbox{6.5cm}{$\phi_y = \eta_t + \phi_x \eta_x$,} y=1+\eta, \label{HW3} \\
& \parbox{6.5cm}{$\phi_t + \frac{1}{2}(\phi_x^2+\phi_y^2)+\eta+\gamma H(\eta)=0$,} y=1+\eta \label{HW4}
\end{align}
with
$$H(\eta) = \frac{1}{(1+\eta_x^2)^{1/2}}\left[\frac{1}{(1+\eta_x^2)^{1/2}}\left(\frac{\eta_{xx}}{(1+\eta_x^2)^{3/2}}\right)_x\right]_x
+\frac{1}{2}\left(\frac{\eta_{xx}}{(1+\eta_x^2)^{3/2}}\right)^3$$
(see Guyenne and Parau \cite{GuyenneParau12}).
Here we have introduced dimensionless variables, choosing $h$ as
length scale and $(h/g)^{1/2}$ as time scale; the parameter $\gamma$ is defined by the formula
$\gamma={\mathcal D}/(\rho g h^4)$, where ${\mathcal D}$, $\rho$ and $g$ are respectively the coefficient of flexural rigidity for the ice sheet,
the density of the fluid and the acceleration due to gravity. \emph{Solitary hydroelastic waves} are non-trivial solutions of these equations of
the form $\eta(x,t)=\eta(x+\nu t)$, $\phi(x,y,t)=\phi(x+\nu t,y)$ with $\eta(x+\nu t) \rightarrow 0$ as $x +\nu t \rightarrow \pm \infty$.

Equations \eqn{HW1}--\eqn{HW4} admit the conserved quantities
$${\mathcal E}(\eta,\Phi) = \frac{1}{2}\int_{-\infty}^\infty \left( \Phi G(\eta) \Phi + \eta^2 + \gamma \frac{\eta_{xx}^2}{(1+\eta_x^2)^{5/2}}\right)\dx,
\qquad
{\mathcal I}(\eta,\Phi) = \int_{-\infty}^\infty \eta_x \Phi \dx
$$
(`energy' and `impulse') associated with translation invariance in $t$ and $x$; the \emph{Dirichlet-Neumann operator} $G(\eta)$ is defined by
$G(\eta)\Phi=(1+\eta_x^2)^{1/2}\phi_n|_{y=1+\eta}$, in which $\phi$ is the harmonic function in ${0<y<1+\eta}$
with $\phi_y|_{y=0}=0$ and $\phi|_{y=1+\eta}=\Phi$.
A hydroelastic solitary wave corresponds to a critical point of the energy under the constraint of fixed impulse
(the potential $\phi$ is recovered from $\Phi$ by
solving the above boundary-value problem)
and therefore a critical point of the functional
${\mathcal E}-\nu {\mathcal I}$, where the Lagrange multiplier $\nu$ gives the wave speed.
Proposition \ref{Analytic} (see Groves \& Wahl\'{e}n \cite[Theorem 2.14(i)]{GrovesWahlen15})
confirms in particular that ${\mathcal E}$, ${\mathcal I}$ are analytic functions
$U \times H_\star^{1/2}({\mathbb R}) \rightarrow {\mathbb R}$,
where $U = B_M(0)$ is a neighbourhood of the origin in $H^2({\mathbb R})$ chosen so that
$U \subseteq W:=\{ \eta \in W^{1,\infty}({\mathbb R}): 1+ \inf_{x \in {\mathbb R}} \eta(x) > h_0\}$
for a fixed $h_0 \in (0,1)$,
and $H_\star^{1/2}({\mathbb R})$,  $H_\star^{-1/2}({\mathbb R})$ are the completions of
${\mathcal S}({\mathbb R})$,
$\overline{\mathcal S}({\mathbb R})=\{\eta \in {\mathcal S}({\mathbb R}): \int_{-\infty}^\infty \eta(x) \dx=0\}$ with respect to the norms $\|\eta\|_{\star,1/2}:=(\int_{-\infty}^\infty (1+k^2)^{-1/2}k^2|\hat{\eta}|^2\dk)^{1/2}$,
$\|\eta\|_{\star,-1/2}:=(\int_{-\infty}^\infty (1+k^2)^{1/2}k^{-2}|\hat{\eta}|^2\dk)^{1/2}$.
\begin{proposition} \label{Analytic}
The mapping $W \rightarrow \mathrm{GL}(H_\star^{1/2}({\mathbb R}),H_\star^{-1/2}({\mathbb R}))$ given by
$\eta \mapsto (\Phi \mapsto G(\eta)\Phi)$ is analytic.
\end{proposition}

Restricting to small-amplitude waves, we seek minimisers of ${\mathcal E}$ subject to the constraint ${\mathcal I}=2\mu$,
where $\mu$ is a small positive number, and establish the following theorem.

\begin{theorem} \label{Main theorem}

The following statements hold for each sufficiently large value of $\gamma$ (see Remark \ref{Sufficiently}).
\begin{enumerate}
\item
The set $D_\mu$ of minimisers of ${\mathcal E}$ over $S_\mu=\{(\eta,\Phi) \in U \times H_\star^{1/2}({\mathbb R}): {\mathcal I}(\eta,\Phi)=2\mu\}$
is non-empty and lies in $H^4({\mathbb R}) \times H_\star^{1/2}({\mathbb R})$. Furthermore, the estimate $\|\eta\|_2 \lesssim \mu^{1/2}$ holds uniformly over $D_\mu$.
\item
Suppose that $\{(\eta_n,\Phi_n)\}$ is a minimising sequence for ${\mathcal E}$.
There exists a sequence $\{x_n\} \subseteq {\mathbb R}$ with the property that a subsequence of
$\{(\eta_n(x_n + \cdot),\Phi_n(x_n + \cdot))\}$ converges in $H^2({\mathbb R})\times H_\star^{1/2}({\mathbb R})$ to a function in $D_\mu$.
\end{enumerate}
\end{theorem}

\begin{remark} (`conditional energetic stability of the set of minimisers')
Suppose that $(\eta,\Phi): [0,T] \rightarrow U \times H_\star^{1/2}({\mathbb R})$
is a solution to \eqn{HW1}--\eqn{HW4} in the sense that
${\mathcal E}(\eta(t),\Phi(t)) = {\mathcal E}(\eta(0),\Phi(0))$, ${\mathcal I}(\eta(t),\Phi(t))={\mathcal I}(\eta(0),\Phi(0))$ for all $t \in [0,T]$
(see Ambrose and Siegel \cite{AmbroseSiegel14} for a discussion of the initial-value problem).
It follows from Theorem \ref{Main theorem} that
for each
$\varepsilon>0$ there exists $\delta>0$ such that $\mathrm{dist}\,((\eta(0),\Phi(0)), D_\mu) < \delta$ implies
$\mathrm{dist}\, ((\eta(t),\Phi(t)), D_\mu)<\varepsilon$ for $t\in[0,T]$, where `dist' denotes the distance in $H^2({\mathbb R})\times H_\star^{1/2}({\mathbb R})$.
\end{remark}

\subsection{Heuristics}

\begin{figure}
\centering
\includegraphics[width=5cm]{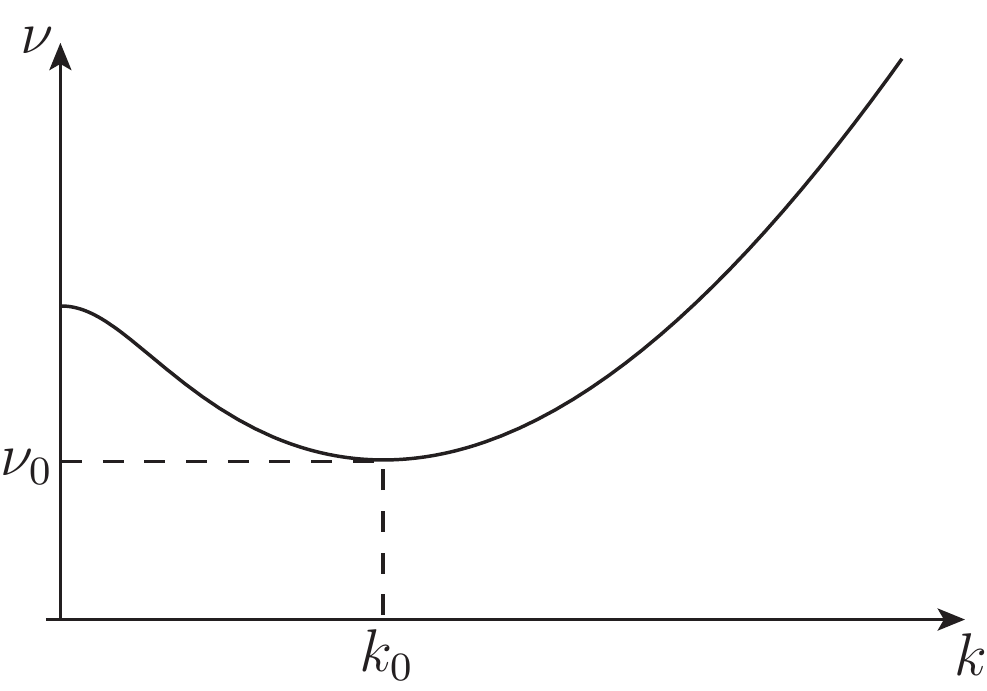}
\hspace{2.5cm}
\includegraphics[width=4cm]{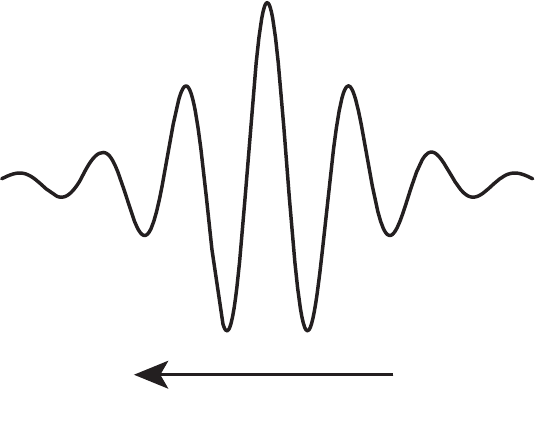}

{\it
\caption{(a) Dispersion relation for linear hydroelastic waves. (b) Small-amplitude envelope solitary waves
with speed $\nu = \nu_0 +2(\nu_0f(k_0))^{-1}\mu^2 \nu_\mathrm{NLS}$ (where $\nu_\mathrm{NLS}<0$)
predicted by nonlinear Schr\"{o}dinger theory.}
\label{Waves}}
\end{figure}

The existence of small-amplitude solitary waves is predicted by studying the dispersion relation for the linearised version of
\eqn{HW1}--\eqn{HW4}. Linear waves of the form $\eta(x,t)=\cos k(x+\nu t)$ exist whenever $\nu=\nu(k)$, where
$\nu(k)^2 = (1+\gamma k^4)/f(k)$, $f(k):=|k|\coth|k|$.
The function $k \mapsto \nu(k)$, $k \geq 0$ has a unique global minimum $\nu_0=\nu(k_0)$ with $k_0>0$ (see Figure 
\ref{Waves}(a)). Note also that $g(k):=1+\gamma k^4-\nu_0^2 f(k) \geq 0$ with equality precisely when
$k=\pm k_0$, and solving the equation $g^\prime(k_0)=0$ yields the relationships $\nu_0^2=4(4f(k_0)-k_0f^\prime(k_0))\big)^{-1}$
and $\gamma=\gamma_0(k_0)$, where
$\gamma_0(k_0)=f^\prime(k_0)\big(k_0^3(4f(k_0)-k_0f^\prime(k_0))\big)^{-1}$,
so that $\gamma_0$ is a strictly monotone decreasing function of $k_0$ with $\lim_{k_0 \rightarrow 0}\gamma(k_0) = \infty$ and
$\lim_{k_0 \rightarrow \infty} \gamma(k_0) =0$.

Bifurcations of nonlinear solitary waves are expected whenever the linear group and phase speeds are equal, so that
$\nu^\prime(k)=0$ (see Dias and Kharif \cite[\S 3]{DiasKharif99}). We therefore expect the existence of small-amplitude
solitary waves with speed near $\nu_0$; the waves bifurcate from a linear periodic wave train with frequency
$k_0\nu_0$ (see Figure \ref{Waves}(b)). The appropriate model equation for this type of solution is the cubic nonlinear
Schr\"{o}dinger equation
\begin{equation}
2\i A_T-\tfrac{1}{4}g^{\prime\prime}(k_0)A_{XX}+\tfrac{3}{2}\left(\tfrac{1}{2}A_3+A_4\right)|A|^2A=0,
\label{NLS}
\end{equation}
in which
$$\eta(x,t)=\tfrac{1}{2}\mu(A(X,T)\e^{\i k_0(x+\nu_0t)} + \mathrm{c.c.}) + O(\mu^2),
\qquad X=\mu(x+\nu_0 t),\ T=2k_0(v_0f(k_0))^{-1}\mu^2t$$
and the abbreviation `c.c.' denotes the complex conjugate of the preceding quantity;
the values of the constants $A_3$ and $A_4$ are
$A_3 = - \tfrac{1}{3}g(2k_0)^{-1}(A_3^1)^2 - \tfrac{2}{3}g(0)^{-1}(A_3^2)^2$ and $A_4=A_4^1-\nu_0^2A_4^2$, where
$$
A_3^1= \nu_0^2 f(2k_0)f(k_0)+\tfrac{1}{2}\nu_0^2f(k_0)^2-\tfrac{3}{2}\nu_0^2k_0^2, \qquad
A_3^2 = \nu_0^2f(k_0)+\tfrac{1}{2}\nu_0^2f(k_0)^2-\tfrac{1}{2}\nu_0^2k_0^2,$$
$$
A_4^1 = -\tfrac{5}{12}\gamma_0 k_0^6, \qquad
A_4^2 = \tfrac{1}{6}f(k_0)^2(f(2k_0)+2) - \tfrac{1}{2}k_0^2f(k_0)
$$
(see Milewski and Wang \cite[\S2]{MilewskiWang13} for a derivation of equation \eqn{NLS}
in the present context).
Note that $k_0>0$; the case $k_0=0$, which is associated with the Korteweg-de Vries scaling limit, does not arise
here.

At this level of approximation, a solution to equation \eqn{NLS} of the form
$A(X,T) = \e^{\i\nu_\mathrm{NLS}T}\zeta(X)$ with $\zeta(X) \rightarrow 0$ as $X \rightarrow \pm \infty$,
so that $\zeta$ is a homoclinic solution of the ordinary differential equation
\begin{equation}
-\tfrac{1}{4}g^{\prime\prime}(k_0)\zeta_{xx}-2\nu_\mathrm{NLS}\zeta+\tfrac{3}{2}\left(
\tfrac{1}{2}A_3+A_4\right)|\zeta|^2\zeta=0 \label{NLS homos}
\end{equation}
with $\nu_\mathrm{NLS}= -\tfrac{9}{8}\alpha_\mathrm{NLS}^2\,g^{\prime\prime}(k_0)^{-1}\left(\tfrac{1}{2}A_3+A_4\right)^2$
and $\alpha_\mathrm{NLS}=2(\nu_0f(k_0))^{-1}$,
corresponds to a solitary wave with speed $\nu = \nu_0 +2(\nu_0f(k_0))^{-1}\mu^2 \nu_\mathrm{NLS}$.

\begin{proposition}
Suppose that $\tfrac{1}{2}A_3 + A_4<0$.
The set of complex-valued homoclinic solutions to the ordinary differential equation
\eqn{NLS homos} is
$D_\mathrm{NLS} = \{\e^{\i\omega}\zeta_\mathrm{NLS}(\cdot + y)\colon \omega \in [0,2\pi), y \in {\mathbb R}\}$, where
$$\zeta_\mathrm{NLS}(x) = \alpha_\mathrm{NLS}\left(-3g^{\prime\prime}(k_0)^{-1}\left(\tfrac{1}{2}A_3+A_4\right)\right)^{\frac{1}{2}}
\sech\left(-3\alpha_\mathrm{NLS}\,g^{\prime\prime}(k_0)^{-1}\left(\tfrac{1}{2}A_3+A_4\right)x\right).$$
\end{proposition}
\begin{remark} \label{Sufficiently}
Since $A_3 < 0$ and $\lim_{k_0 \rightarrow 0}A_4 = -\frac{1}{2}$, so that $A_4 < 0$ for sufficiently small values of $k_0$,
we find that $\tfrac{1}{2}A_3 + A_4<0$ for sufficiently small values of $k_0$, or
equivalently for sufficiently large values of $\gamma$ (corresponding to sufficiently shallow water
in physical variables). Numerics indicate that $\tfrac{1}{2}A_3 + A_4<0$ for $k_0<177.33$, or equivalently $\gamma>3.37 \times 10^{-10}$.
\end{remark}

Our second theorem confirms the heuristic argument given above.
\begin{theorem} \label{Convergence}
Suppose that $\frac{1}{2}A_3 + A_4<0$. The set $D_\mu$ of minimisers of ${\mathcal E}$ over $S_\mu$ satisfies
$$\sup_{(\eta,\Phi) \in D_\mu} \inf_{\omega \in [0,2\pi], x \in {\mathbb R}} \|\zeta_\eta-e^{\i\omega}\zeta_\mathrm{NLS}(\cdot+x)\|_1 \rightarrow 0$$
as $\mu \downarrow 0$, where we write
$\eta_1^+(x) = \tfrac{1}{2}\mu \zeta_\eta(\mu x)\e^{\i k_0 x}$
and $\eta_1^+={\mathcal F}^{-1}[\chi_{[k_0-\delta_0,k_0+\delta_0]}\hat{\eta}]$ with $\delta_0 \in (0,\tfrac{1}{3}k_0)$.
Furthermore, the speed $\nu_\mu$ of the corresponding solitary
wave satisfies
$\nu_\mu = \nu_0 + 2(\nu_0f(k_0))^{-1}\nu_\mathrm{NLS}\mu^2 + o(\mu^2)$
uniformly over $(\eta,\Phi) \in D_\mu$.
\end{theorem}

\section{The constrained minimisation problem} \label{Minimisation problem}

We tackle the constrained minimisation problem in two steps.
(i) Fix $\eta \neq 0$ and minimise ${\mathcal E}(\eta,\cdot)$ over
$T_\mu=\{\Phi \in H^{1/2}_\star({\mathbb R}): {\mathcal I}(\eta,\Phi)=2\mu\}$.
This problem (of minimising a quadratic functional over a linear manifold)
admits a unique global minimiser $\Phi_\eta$.
(ii) Minimise ${\mathcal J}_\mu(\eta):={\mathcal E}(\eta,\Phi_\eta)$ over $\eta \in U\setminus\{0\}$.
Because $\Phi_\eta$ minimises
${\mathcal E}(\eta,\cdot)$ over $T_\mu$ there exists a Lagrange multiplier $\nu_\eta$ such that
$G(\eta)\Phi_\eta = \nu_\eta \eta_x$,
and straightforward calculations show that
$\Phi_\eta = \nu_\eta G(\eta)^{-1}\eta_x$, $\nu_\eta = \mu/{\mathcal L}(\eta)$
and
$$
{\mathcal J}_\mu(\eta)={\mathcal K}(\eta)+\frac{\mu^2}{{\mathcal L}(\eta)},
$$
where
$$
{\mathcal K}(\eta) = \frac{1}{2}\int_{-\infty}^\infty
\left(\eta^2 + \frac{\gamma\eta_{xx}^2}{(1+\eta_x^2)^{5/2}}\right) \dx, \qquad
{\mathcal L}(\eta) = \frac{1}{2}\int_{-\infty}^\infty \eta_x G(\eta)^{-1}\eta_x \dx.
$$
This computation also shows that the dimensionless speed of the solitary wave corresponding to
a constrained minimiser of ${\mathcal E}$ over $S_\mu$ is $\mu/{\mathcal L}(\eta)$.

A similar minimisation problem arises in the study of irrotational solitary water waves with weak surface tension
(see Groves and Wahl\'{e}n \cite{GrovesWahlen15}, taking $\omega=0$ and $\beta<\beta_\mathrm{c}$); in that case
${\mathcal K}(\eta)$ is replaced by $\tilde{\mathcal K}(\eta) = \int_{-\infty}^\infty \big(\tfrac{1}{2}\eta^2 + \beta((1+\eta_x^2)^{1/2}-1)\big)\dx$. In this note we describe the modifications necessary to apply the theory of Groves and
Wahl\'{e}n to the hydroelastic problem. The presence of the second-order derivative necessitates on the one
hand non-trivial modifications because the $L^2({\mathbb R})$-gradient ${\mathcal K}^\prime(\eta)$ is not defined on
the whole of $U$, but leads on the other hand to a more satisfactory final result (compare Theorem
\ref{Main theorem} with Theorem 1.5 of Groves \& Wahl\'{e}n).

Lemmata \ref{Basic properties of K} and \ref{Basic properties of L},  in which we write
$W^s:=W \cap H^s({\mathbb R})$, state some basic properties of the functionals ${\mathcal K}$ and ${\mathcal L}$
(see Groves and Wahl\'{e}n \cite{GrovesWahlen15} for the proof of the latter), while Proposition \ref{Technical result} is a useful
`weak-strong' argument. Note that the `linear' estimates for ${\mathcal K}_\mathrm{nl}(\eta)$ and ${\mathcal L}_\mathrm{nl}(\eta)$
are used only to bound the $W^{1,\infty}({\mathbb R})$ norm of a minimising sequence for ${\mathcal J}$ over $U \setminus \{0\}$
away from zero (see the discussion at the beginning of Section 3).

\begin{lemma} \quad \label{Basic properties of K}
\begin{enumerate}
\item
The functional ${\mathcal K}: H^2({\mathbb R}) \rightarrow {\mathbb R}$ is analytic and satisfies ${\mathcal K}(0)=0$.
\item
There exists a constant $D>0$ such that ${\mathcal K}(\eta) \geq D^{-1} \|\eta\|_2^2$ for all $\eta \in U$.
\item
The $L^2({\mathbb R})$-gradient ${\mathcal K}^\prime(\eta)$ exists for each $\eta \in H^4({\mathbb R})$ and is
given by the formula
$${\mathcal K}^\prime(\eta) = \eta + \gamma \left[\frac{\eta_{xx}}{(1+\eta_x^2)^{5/2}}\right]_{xx}
+ \frac{5}{2}\gamma \left[\frac{\eta_x\eta_{xx}^2}{(1+\eta_x^2)^{7/2}}\right]_x.$$
This formula defines an analytic function ${\mathcal K}^\prime: H^2({\mathbb R}) \rightarrow H^{-2}({\mathbb R})$ which
satisfies ${\mathcal K}^\prime(0)=0$.
\item
The estimates
$|{\mathcal K}_4(\eta)|
\lesssim \|\eta\|_2^2 \|\eta\|_{1,\infty}^2$,
$|{\mathcal K}_\mathrm{r}(\eta)|
\lesssim \|\eta\|_2^3 \|\eta\|_{1,\infty}^2$,
$|{\mathcal K}_\mathrm{nl}(\eta)| \lesssim \|\eta\|_{1,\infty}$
hold for all $\eta \in U$, where ${\mathcal K}_n(\eta)=\frac{1}{n!}\mathrm{d}^n{\mathcal K}[0](\{\eta\}^n)$,
${\mathcal K}_\mathrm{r}(\eta)=\sum_{n=5}^\infty {\mathcal K}_n(\eta)$ and ${\mathcal K}_\mathrm{nl}(\eta)={\mathcal K}(\eta)-{\mathcal K}_2(\eta)$.
\item
The estimates
$$\|{\mathcal F}^{-1}[(1-\chi_S(k))g(k)^{-1/2}{\mathcal F}[{\mathcal K}_4^\prime(\eta)]]\|_0 \lesssim \|\eta\|_2 (\|\eta\|_{1,\infty} + \|\eta_{xx}+k_0^2\eta\|_0)^2,$$
$$\|{\mathcal F}^{-1}[(1-\chi_S(k))g(k)^{-1/2}{\mathcal F}[{\mathcal K}_\mathrm{r}^\prime(\eta)]]\|_0,
\ |\langle {\mathcal K}_4^\prime(\eta),\eta \rangle_0|,\ |\langle {\mathcal K}_\mathrm{r}^\prime(\eta),\eta \rangle_0|
 \lesssim\|\eta\|_2^2 (\|\eta\|_{1,\infty} + \|\eta_{xx}+k_0^2\eta\|_0)^2$$
hold for all $\eta \in H^2({\mathbb R})$, where $S=[-k_0-\delta_0,-k_0+\delta_0] \cup [k_0-\delta_0,k_0+\delta_0]$ and
$\delta_0 \in (0,\tfrac{1}{3}k_0)$.
\end{enumerate}
\end{lemma}
{\bf Proof.} Assertions (i)--(iv) follow by straightforward estimates. Turning to (v), note that 
$$
{\mathcal K}_4^\prime(\eta) = \tfrac{5}{2}\gamma\big((\eta_x\eta_{xx}^2)_x+(\eta_x^2\eta_{xx})_{xx}\big)
= \tfrac{5}{2}\gamma\left(\big(\eta_x(\eta_{xx}+k_0^2\eta)^2-2k_0^2\eta_x\eta(\eta_{xx}+k_0^2\eta)+k_0^4\eta_x\eta^2\big)_x
+\big(\eta_x^2\eta_{xx})_{xx}\right)
$$
so that
\begin{align*}
\|{\mathcal F}^{-1}[(1-\chi_S(k))g(k)^{-1/2}{\mathcal F}[{\mathcal K}_4^\prime(\eta)]]\|_0
& \lesssim \|\eta_x(\eta_{xx}+k_0^2\eta)^2\|_{-1}+\|\eta_x\eta(\eta_{xx}+k_0^2\eta)\|_0+\|\eta_x\eta^2\|_0
+\|\eta_x^2\eta_{xx}\|_0 \\
& \lesssim \|\eta_x(\eta_{xx}+k_0^2\eta)\|_0\|\eta_{xx}+k_0^2\eta\|_0+\|\eta\|_2 \|\eta\|_{1,\infty}^2 \\
& \lesssim \|\eta\|_2 (\|\eta\|_{1,\infty} + \|\eta_{xx}+k_0^2\eta\|_0)^2,
\end{align*}
where we have used the inequalities $(1-\chi_S(k))g(k)^{-1/2} \lesssim (1+|k|^2)^{-1}$ and $\|u_1u_2\|_{-1} \lesssim \|u_1\|_0\|u_2\|_0$
(see H\"{o}rmander \cite[Theorem 8.3.1]{Hoermander}); the remaining estimates are obtained in a similar fashion.\qed\pagebreak

\begin{lemma} \quad \label{Basic properties of L}
\begin{enumerate}
\item
Suppose $s>0$. The functional ${\mathcal L}: W^{s+3/2} \rightarrow {\mathbb R}$
is analytic and satisfies ${\mathcal L}(0)=0$.
\item
The estimates $\|\eta\|_{1/2}^2 \lesssim {\mathcal L}(\eta), {\mathcal L}_2(\eta) \lesssim \|\eta\|_{1/2}^2$,
where ${\mathcal L}_2(\eta)=\frac{1}{2!}\mathrm{d}^2{\mathcal L}[0](\{\eta\}^2)$, hold for all $\eta \in U$.
\item
Suppose $s>0$.
The $L^2({\mathbb R})$-gradient ${\mathcal L}^\prime(\eta)$ exists for each $\eta \in W^{s+3/2}$ and defines
an analytic function
${\mathcal L}^\prime: W^{s+3/2} \rightarrow H^{s+1/2}({\mathbb R})$ which satisfies
${\mathcal L}^\prime(0)=0$.
\item
Suppose that $\{M_n^{(1)}\}$, $\{M_n^{(2)}\} \subseteq {\mathbb R}$ and $\{\eta_n^{(1)}\}$, $\{\eta_n^{(2)}\} \subseteq U$ are
sequences with $M_n^{(1)}$, $M_n^{(2)} \rightarrow \infty$, $M_n^{(1)}/M_n^{(2)} \rightarrow 0$,
$\{\eta_n^{(1)}+\eta_n^{(2)}\} \subseteq U$ and
$\mathrm{supp}\, \eta_n^{(1)} \subseteq (-2M_n^{(1)},2M_n^{(1)})$,
$\mathrm{supp}\, \eta_n^{(2)} \subseteq {\mathbb R} \setminus (-M_n^{(2)},M_n^{(2)})$. The functional ${\mathcal L}$ has
the `pseudolocal' properties
$${\mathcal L}(\eta_n^{(1)}+\eta_n^{(2)})-{\mathcal L}(\eta_n^{(1)})-{\mathcal L}(\eta_n^{(2)}) \rightarrow 0, \qquad
\|{\mathcal L}^\prime(\eta_n^{(1)}+\eta_n^{(2)})-{\mathcal L}^\prime(\eta_n^{(1)})-{\mathcal L}^\prime(\eta_n^{(2)})\|_0 \rightarrow 0$$
and $\langle {\mathcal L}^\prime(\eta_n^{(2)}), \phi \rangle_0 \rightarrow 0$ for each $\phi \in C_0^\infty({\mathbb R})$.
\item
The estimates
$$
|{\mathcal L}_3(\eta)|
\lesssim \|\eta\|_2^2 (\|\eta\|_{1,\infty} + \|\eta_{xx}+k_0^2\eta\|_0), \qquad
|{\mathcal L}_4(\eta)|
\lesssim \|\eta\|_2^2 (\|\eta\|_{1,\infty} + \|\eta_{xx}+k_0^2\eta\|_0)^2,
$$
$$
|{\mathcal L}_\mathrm{r}(\eta)|
\lesssim \|\eta\|_2^3 (\|\eta\|_{1,\infty} + \|\eta_{xx}+k_0^2\eta\|_0)^2, \qquad
|{\mathcal L}_\mathrm{nl}(\eta)|
\lesssim \|\eta\|_{1,\infty},
$$
where ${\mathcal L}_n(\eta)=\frac{1}{n!}\mathrm{d}^n{\mathcal L}[0](\{\eta\}^n)$,
${\mathcal L}_\mathrm{r}(\eta)=\sum_{n=5}^\infty {\mathcal L}_n(\eta)$ and
${\mathcal L}_\mathrm{nl}(\eta)={\mathcal L}(\eta)-{\mathcal L}_2(\eta)$, and
\begin{align*}
\|{\mathcal L}_3^\prime(\eta)\|_0
& \lesssim \|\eta\|_2 (\|\eta\|_{1,\infty} + \|\eta_{xx}+k_0^2\eta\|_0 + \|K^0 \eta\|_\infty), \\
\|{\mathcal L}_4^\prime(\eta)\|_0
& \lesssim \|\eta\|_2 (\|\eta\|_{1,\infty} + \|\eta_{xx}+k_0^2\eta\|_0 + \|K^0 \eta\|_\infty)^2, \\
\|{\mathcal L}_\mathrm{r}^\prime(\eta)\|_0
& \lesssim \|\eta\|_2^2  (\|\eta\|_{1,\infty} + \|\eta_{xx}+k_0^2\eta\|_0)^2,
\end{align*}
where $K_0\eta:={\mathcal F}^{-1}[f(k)\hat{\eta}]$, hold for all $\eta \in U$.
\end{enumerate}
\end{lemma}

\begin{proposition} \label{Technical result}
Suppose that $\{\eta_n\} \subseteq U$ and $\eta \in U$ have the properties that
$\eta_n \rightharpoonup \eta$ in $H^2({\mathbb R})$ and
$\eta_n \rightarrow \eta$ in $L^2({\mathbb R})$ (and hence in $H^s({\mathbb R})$ for all $s \in [0,2)$). The inequality
${\mathcal K}(\eta) \leq \lim_{n \rightarrow \infty} {\mathcal K}(\eta_n)$ holds whenever $\{{\mathcal K}(\eta_n)\}$ is convergent,
and equality implies that $\eta_n \rightarrow \eta$ in $H^2({\mathbb R})$.
\end{proposition}
{\bf Proof.} Note that
$(1+\eta_{nx}^2)^{-5/4}\eta_{nxx} \rightharpoonup (1+\eta_x^2)^{-5/4}\eta_{xx}$ in $L^2({\mathbb R})$, and it follows from
the weak lower semicontinuity of $\|\cdot\|_0^2$ (and $\eta_n \rightarrow \eta$ in $L^2({\mathbb R})$) that 
${\mathcal K}(\eta) \leq \lim_{n \rightarrow \infty} {\mathcal K}(\eta_n)$. Moreover, ${\mathcal K}(\eta_n) \rightarrow {\mathcal K}(\eta)$
implies that
$\|(1+\eta_{nx}^2)^{-5/4}\eta_{nxx}\|_0 \rightarrow \|(1+\eta_x^2)^{-5/4}\eta_{xx}\|_0$, so that
$(1+\eta_{nx}^2)^{-5/4}\eta_{nxx} \rightarrow (1+\eta_x^2)^{-5/4}\eta_{xx}$ in $L^2({\mathbb R})$ and hence
$\eta_{nxx} \rightarrow \eta_{xx}$ in $L^2({\mathbb R})$.\qed

Next we establish some basic properties of ${\mathcal J}_\mu$. The following proposition (cf.\ Groves and Wahl\'{e}n \cite[Appendix A.2]{GrovesWahlen15})
shows in particular that $c_\mu:=\inf_{\eta \in U \setminus \{0\}} {\mathcal J}_\mu(\eta)<2\nu_0\mu$, while Lemma \ref{Regularity} shows that its critical points have additional regularity.

\begin{proposition} \label{Test function}
The continuous mapping $\alpha \mapsto \nu_0 {\mathcal L}(\eta^\star_\alpha)$, where
$$\eta^\star_\alpha(x) = \alpha \zeta_\mathrm{NLS}(\alpha x) \cos k_0 x - \tfrac{1}{2}\alpha^2g(2k_0)^{-1}A_3^1 \zeta_\mathrm{NLS}(\alpha x)^2 \cos 2k_0 x - \tfrac{1}{2}\alpha^2g(0)^{-1}A_3^2 \zeta_\mathrm{NLS} (\alpha x)^2,$$
is invertible, and its (continuous) inverse $\mu \mapsto \alpha(\mu)$ satisfies
${\mathcal J}_\mu(\eta^\star_{\alpha(\mu)}) = 2\nu_0\mu + c_\mathrm{NLS}\mu^3 + o(\mu^3)$, where
$c_\mathrm{NLS}=-\frac{3}{4}\alpha_\mathrm{NLS}^3g^{\prime\prime}(k_0)^{-1}(\frac{1}{2}A_3+A_4)^2$.
\end{proposition}

\begin{remark}
Each $\eta \in U \setminus \{0\}$ satisfies
$${\mathcal K}_2(\eta) + \frac{\mu^2}{{\mathcal L}_2(\eta)} 
={\mathcal K}_2(\eta) -\nu_0^2 {\mathcal L}_2(\eta)+\frac{(\mu-\nu_0{\mathcal L}_2(\eta))^2}{{\mathcal L}_2(\eta)}+2\nu_0\mu
\geq \frac{1}{2}\int_{-\infty}^\infty g(k) |\hat{\eta}|^2 \dk + 2\nu_0\mu \geq 2 \nu_0\mu.
$$

\end{remark}

\begin{lemma} \label{Regularity}
Any critical point $\eta \in U \setminus \{0\}$ of ${\mathcal J}_\mu$ belongs to $H^4({\mathbb R})$.
\end{lemma}
{\bf Proof.} 
Write $u=(1+\eta_x^2)^{-5/2}\eta_{xx}$, so that $\eta_x(1+\eta_x^2)^{3/2}u^2 \in L^1({\mathbb R}) \subseteq H^{-3/4}({\mathbb R})$, and observe that
\begin{equation}
\gamma u_{xx} = \frac{\mu}{{\mathcal L}(\eta)^2}{\mathcal L}^\prime(\eta) - \eta - \tfrac{5}{2}\gamma\big(\eta_x(1+\eta_x^2)^{3/2}u^2\big)_x \label{Bootstrap}
\end{equation}
in the sense of distributions since $\eta$ is a critical point of ${\mathcal J}_\mu$.
It follows from \eqn{Bootstrap} and the fact that ${\mathcal L}^\prime(\eta) \in L^2({\mathbb R})$ that $\gamma u_{xx} \in H^{-7/4}({\mathbb R})$, that
is $u \in H^{1/4}({\mathbb R})$. We conclude that $u^2 \in L^2({\mathbb R})$ (see H\"{o}rmander \cite[Theorem 8.3.1]{Hoermander}), so that
$\eta_x(1+\eta_x^2)^{3/2}u^2 \in L^2({\mathbb R})$ and hence $\gamma u_{xx} \in H^{-1}({\mathbb R})$, that is $u \in H^1({\mathbb R})$.

Observing that $\eta_x(1+\eta_x^2)^{3/2}u^2 \in H^1({\mathbb R})$, one finds from \eqn{Bootstrap} that
$\gamma u_{xx} \in L^2({\mathbb R})$, $u \in H^2({\mathbb R})$ and finally $\eta \in H^4({\mathbb R})$
(because $\eta_{xx}=(1+\eta_x^2)^{5/2}u$).\qed\pagebreak

Theorem \ref{Main theorem} is a consequence of the following result (cf.\ Groves \& Wahl\'{e}n
\cite[Theorem 5.2]{GrovesWahlen15}).

\begin{theorem} \label{Implies main theorem}
Suppose that $\frac{1}{2}A_3 + A_4<0$.
\begin{enumerate}
\item
The set $B_\mu$ of minimisers of ${\mathcal J}_\mu$ over $U \setminus \{0\}$ is nonempty and lies in $H^4({\mathbb R})$.
Moreover, each $\eta \in B_\mu$ satisfies $\|\eta\|_2^2 \leq 2D\nu_0\mu$.
\item
Suppose that $\{\eta_n\}$ is a minimising sequence for ${\mathcal J}_\mu$ over $U \setminus \{0\}$. There exists a sequence
$\{x_n\} \subseteq {\mathbb R}$ with the property that there exists a subsequence of $\{\eta_n(x_n+\cdot)\}$ which converges
in $H^2({\mathbb R})$ to a function $\eta \in B_\mu$.
\end{enumerate}
\end{theorem}

Any function $\eta \in U$ with ${\mathcal J}_\mu(\eta) < 2\nu_0\mu$ satisfies
$\|\eta\|_2^2 < 2D\nu_0\mu$, ${\mathcal L}(\eta) > \mu/(2\nu_0)$ and
${\mathcal L}_2(\eta) \gtrsim \mu$ (see Lemmata \ref{Basic properties of K}(ii) and \ref{Basic properties of L}(ii)).
These properties are enjoyed in particular by a minimising sequence $\{\eta_n\}$ for ${\mathcal J}_\mu$ over $U \setminus \{0\}$,
which also satisfies ${\mathcal M}_\mu(\eta_n) \lesssim - \mu^3$, where
${\mathcal M}_\mu(\eta) = {\mathcal J}_\mu(\eta) - {\mathcal K}_2(\eta) - \mu^2/{\mathcal L}_2(\eta)$
(Proposition \ref{Test function}), and hence $\|\eta_n\|_{1,\infty} \gtrsim \mu^3$ (because $|{\mathcal K}_\mathrm{nl}(\eta_n)|$,
$|{\mathcal L}_\mathrm{nl}(\eta_n)| \lesssim \|\eta_n\|_{1,\infty}$). Furthermore, we may without loss of generality assume that
$\{\eta_n\}$ is a Palais-Smale sequence, so that $\mathrm{d}{\mathcal J}_\mu[\eta_n] \rightarrow 0$
in $(H^2({\mathbb R}))^\ast$, and the calculation
$$\|{\mathcal J}^\prime(\eta_n)\|_{-2} = \sup \{ \langle {\mathcal J}^\prime(\eta_n), \phi \rangle_0: \phi \in H^2({\mathbb R}),
\|\phi\|_2 =1\} = \|\mathrm{d}{\mathcal J}_\mu[\eta_n]\|_{(H^2({\mathbb R}))^\ast}$$
shows that ${\mathcal J}^\prime(\eta_n) \rightarrow 0$ in $H^{-2}({\mathbb R})$.
Theorem \ref{Implies main theorem} is proved by applying the concentration-compactness principle to
the sequence $\{\eta_{nx}^2 + \eta_n^2\} \subseteq L^1({\mathbb R})$ under the additional hypothesis that $c_\mu$ is
a strictly sub-additive function of $\mu$, which is verified in Section \ref{SSA} below.

\emph{`Vanishing'} is excluded since it implies that $\|\eta_n\|_{1,\infty} \rightarrow 0$, which contradicts the estimate
$\|\eta_n\|_{1,\infty} \gtrsim \mu^3$ (see above).

\emph{`Dichotomy'} leads to the existence of sequences
$\{\eta_n^{(1)}\}$, $\{\eta_n^{(2)}\}$ of the kind described in Lemma \ref{Basic properties of L}(iv) 
with
$\lim_{n \rightarrow \infty} \|\eta_n - \eta_n^{(1)}-\eta_n^{(2)}\|_2 = 0$  (up to subsequences and translations), so that in particular
$$\lim_{n \rightarrow \infty} {\mathcal J}_\mu(\eta_n) = \lim_{n \rightarrow \infty} {\mathcal J}_{\mu^{(1)}}(\eta_n^{(1)})+
\lim_{n \rightarrow \infty} {\mathcal J}_{\mu^{(2)}}(\eta_n^{(2)}),$$
where
$\mu^{(j)}\ =\ \mu\lim_{n \rightarrow \infty} {\mathcal L}(\eta_n^{(j)})/\lim_{n \rightarrow \infty} \vphantom{{\mathcal L}^2}{\mathcal L}(\eta_n)$
(so that $\mu^{(1)}+\mu^{(2)}=\mu$). We thus obtain the contradiction
$$c_\mu < c_{\mu^{(1)}}+c_{\mu^{(2)}} \leq  \lim_{n \rightarrow \infty} {\mathcal J}_{\mu^{(1)}}(\eta_n^{(1)})+
\lim_{n \rightarrow \infty} {\mathcal J}_{\mu^{(2)}}(\eta_n^{(2)}) = \lim_{n \rightarrow \infty} {\mathcal J}_\mu(\eta_n) = c_\mu,$$
which excludes `dichotomy'.

\emph{`Concentration'} implies the existence
of $\eta \in U$ with
$\eta_n \rightharpoonup \eta$ in $H^2({\mathbb R})$ and
$\eta_n \rightarrow \eta$ in $L^2({\mathbb R})$ (up to subsequences and translations).
Since ${\mathcal K}(\eta_n) \leq {\mathcal J}_\mu(\eta_n) < 2\nu_0\mu$ the sequence $\{{\mathcal K}(\eta_n)\}$ is bounded
and hence admits a convergent subsequence (still denoted by $\{{\mathcal K}(\eta_n)\}$) which satisfies ${\mathcal K}(\eta)
\leq \lim_{n \rightarrow \infty} {\mathcal K}(\eta_n)$ (Proposition \ref{Technical result}).
Lemma \ref{Basic properties of L}(i) asserts that ${\mathcal L}(\eta_n) \rightarrow {\mathcal L}(\eta)$, so that
${\mathcal J}_\mu(\eta) \leq \lim_{n \rightarrow \infty} {\mathcal J}(\eta_n)=c_\mu$, which therefore holds with equality;
it follows that ${\mathcal K}(\eta_n) \rightarrow {\mathcal K}(\eta)$ and hence $\eta_n \rightarrow \eta$ in $H^2({\mathbb R})$
(Proposition \ref{Technical result}), so that $\eta$ minimises ${\mathcal J}_\mu$ over $U \setminus \{0\}$.

\section{Strict sub-additivity} \label{SSA}

We begin
by deriving sharper estimates for a `near minimiser' of ${\mathcal J}_\mu$ over $U\setminus\{0\}$, that is a function
$\tilde{\eta} \in U\setminus\{0\}$ with $\|{\mathcal J}_\mu^\prime(\tilde{\eta})\|_{-2} \leq \mu^N$ for some $N \in {\mathbb N}$
and ${\mathcal J}_\mu(\tilde{\eta})<2\nu_0\mu$ (and hence $\|\tilde{\eta}\|_2 \lesssim \mu^{1/2}$,
${\mathcal L}(\tilde{\eta})$, ${\mathcal L}_2(\tilde{\eta}) \geq \mu$); these estimates apply in particular to a minimising sequence
$\{\eta_n\}$ for ${\mathcal J}_\mu$ over $U \setminus \{0\}$.

We write the equation ${\mathcal J}^\prime_\mu(\eta)={\mathcal K}^\prime(\eta)
-(\mu/{\mathcal L}(\eta))^2{\mathcal L}^\prime(\eta)$ for $\eta \in U$ in the form
$$g(k)\hat{\eta} = {\mathcal F} \left[{\mathcal J}^\prime_\mu(\eta)-{\mathcal K}_\mathrm{nl}^\prime(\eta)
+\left(\frac{\mu}{{\mathcal L}(\eta)}+\nu_0\right)\!\!\!\left(\frac{\mu}{{\mathcal L}(\eta)}-\nu_0\right){\mathcal L}_2^\prime(\eta)
+\left(\frac{\mu}{{\mathcal L}(\eta)}\right)^2{\mathcal L}_\mathrm{nl}^\prime(\eta)\right]
$$
and decompose it into two coupled equations by defining $\eta_2 \in H^2({\mathbb R})$ by the formula
$$\eta_2 = {\mathcal F}^{-1}\left[\frac{1-\chi_S(k)}{g(k)}{\mathcal F} \left[{\mathcal J}^\prime_\mu(\eta)-{\mathcal K}_\mathrm{nl}^\prime(\eta)
+\left(\frac{\mu}{{\mathcal L}(\eta)}+\nu_0\right)\!\!\!\left(\frac{\mu}{{\mathcal L}(\eta)}-\nu_0\right){\mathcal L}_2^\prime(\eta)
+\left(\frac{\mu}{{\mathcal L}(\eta)}\right)^2{\mathcal L}_\mathrm{nl}^\prime(\eta)\right]\right]
$$
(recall that $(1-\chi_S(k))g(k)^{-1/2} \lesssim (1+|k|^2)^{-1}$)
and $\eta_1 \in H^2({\mathbb R})$ by $\eta_1=\eta-\eta_2$, so that $\mathrm{supp}\, \hat{\eta}_1 \in S$ and
$\chi_S {\mathcal L}_3^\prime(\eta_1)=0$ (see Groves and Wahl\'{e}n \cite[Proposition 4.15]{GrovesWahlen15}). We accordingly write these equations as
$$g(k)\hat{\eta}_1=\chi_S(k){\mathcal F}[{\mathcal R}(\eta)-{\mathcal K}_\mathrm{nl}^\prime(\eta)],
\qquad \eta_3:=\eta_2+H(\eta)={\mathcal F}^{-1}\left[\frac{1-\chi_S(k)}{g(k)}{\mathcal F}[{\mathcal R}(\eta)-{\mathcal K}_\mathrm{nl}^\prime(\eta)]\right],$$
where
$$H(\eta):={\mathcal F}^{-1}\left[\frac{1}{g(k)}
{\mathcal F}\left[-\Big(\frac{\mu}{{\mathcal L}(\eta)}\Big)^2{\mathcal L}_3^\prime(\eta_1)\right]\right],$$
$$
{\mathcal R}(\eta):=
{\mathcal J}^\prime_\mu(\eta)
+\left(\frac{\mu}{{\mathcal L}(\eta)}+\nu_0\right)\!\!\!\left(\frac{\mu}{{\mathcal L}(\eta)}-\nu_0\right){\mathcal L}_2^\prime(\eta)
+\left(\frac{\mu}{{\mathcal L}(\eta)}\right)^2({\mathcal L}_\mathrm{nl}^\prime(\eta)-{\mathcal L}_3^\prime(\eta_1)).
$$

The next step is to study $\eta_1$ using the scaled norm
$$
\nn \eta_1 \nn_\alpha := 
\left(\int _{-\infty}^\infty (1+\mu^{-4\alpha}(|k|-k_0)^4) |\hat{\eta}_1(k)|^2\dk\right)^{1/2}
$$
for $H^2({\mathbb R})$; we choose $\alpha>0$ as large as possible so that $\nn \tilde{\eta}_1 \nn_\alpha \lesssim \mu^{1/2}$.

\begin{lemma}
Each near minimiser $\tilde{\eta}$ of ${\mathcal J}_\mu$ over $U \setminus \{0\}$ satisfies
$\|H(\tilde{\eta})\|_2 \lesssim \mu^{1/2+\alpha/2}\nn \tilde{\eta}_1 \nn_\alpha$,\linebreak
$\|{\mathcal R}(\tilde{\eta})\|_{-2} \lesssim \mu^{1/2+\alpha}\nn \tilde{\eta}_1\nn_\alpha^2 + \mu^N$
and
$\|{\mathcal F}^{-1}[(1-\chi_S(k))g(k)^{-1/2}{\mathcal F}[{\mathcal K}_\mathrm{nl}^\prime(\tilde{\eta})]]\|_0 \lesssim \mu^{1/2+\alpha}\nn \tilde{\eta}_1 \nn_\alpha^2+\mu\|\tilde{\eta}_3\|_2$.
\end{lemma}
{\bf Proof.} The results for $H(\tilde{\eta})$ and ${\mathcal R}(\tilde{\eta})$ were derived by Groves \& Wahl\'{e}n \cite[\S4.3.1]{GrovesWahlen15}, while that
for ${\mathcal K}_\mathrm{nl}^\prime(\tilde{\eta})$ follows from Lemma \ref{Basic properties of K}(v) and 
the estimates $\|\eta_1\|_{1,\infty} \lesssim \mu^{\alpha/2}\nn \eta_1 \nn_\alpha$
and $\|\eta_{1xx}+k_0^2\eta_1\|_0 \leq c \mu^{\alpha}\nn \eta_1\nn_\alpha$
(Groves and Wahl\'{e}n \cite[Proposition 4.1]{GrovesWahlen15}).\qed

Square integrating the equation $g(k)\hat{\eta}_1=\chi(k){\mathcal F}[{\mathcal R}(\eta)-{\mathcal K}_\mathrm{nl}^\prime(\eta)]$,
multiplying by $\mu^{-4\alpha}$ and adding $\|\tilde{\eta}_1\|_0^2 \lesssim \mu$
yields $\nn \tilde{\eta}_1 \nn_\alpha^2 \lesssim \mu^{1-2\alpha} \nn \tilde{\eta}_1 \nn_\alpha^4 + \mu$,
which implies that $\nn \tilde{\eta}_1 \nn_\alpha^2 \lesssim \mu$ for each $\alpha<1$;
it follows that
$\|\tilde{\eta}_3\|_2^2 \lesssim\mu^{3+2\alpha}$ and $\|H(\tilde{\eta})\|_2^2 \lesssim\mu^{2+\alpha}$ for each $\alpha<1$.
These estimates are used to establish the following proposition (see Groves \& Wahl\'{e}n \cite[\S4.3.2]{GrovesWahlen15}).

\begin{proposition}
Suppose that $\tilde{\eta}$ is a near minimiser of ${\mathcal J}_\mu$ over $U \setminus \{0\}$. The estimates
\begin{align*}
{\mathcal M}_{a^2\mu}(a\tilde{\eta}) &= -a^3\nu_0^2{\mathcal L}_3(\tilde{\eta})-a^4\nu_0^2{\mathcal L}_4(\tilde{\eta}) + a^3 o(\mu^3),\\
\langle {\mathcal M}_{a^2\mu}^\prime(a\tilde{\eta}), a\tilde{\eta} \rangle_0 + 4a^2\mu \tilde{{\mathcal M}}_{a^2\mu}(a\tilde{\eta})
&=-3a^3\nu_0^2{\mathcal L}_3(\tilde{\eta})-4a^4\nu_0^2{\mathcal L}_4(\tilde{\eta}) + a^3 o(\mu^3),
\end{align*}
where $\tilde{\mathcal M}_\mu(\eta)=\mu/{\mathcal L}(\eta)-\mu/{\mathcal L}_2(\eta)$,
hold uniformly over $a \in [1,2]$.
\end{proposition}

\begin{lemma}
Each near minimiser $\tilde{\eta}$ of ${\mathcal J}_\mu$ over $U \setminus \{0\}$ satisfies the estimate
$${\mathcal K}_4(\tilde{\eta}) = A_4^1 \int_{-\infty}^\infty \tilde{\eta}_1^4\dx+ o(\mu^3).$$
\end{lemma}
{\bf Proof.} We expand the right-hand side of the formula
$${\mathcal K}_4(\tilde{\eta}) = -\frac{5}{4}\gamma \int_{-\infty}^\infty (\partial_x(\tilde{\eta}_1 + H(\tilde{\eta}) + \tilde{\eta}_3))^2 \partial_x^2((\tilde{\eta}_1 + H(\tilde{\eta}) + \tilde{\eta}_3))^2\dx;$$
terms with zero, one or two occurrences of $\tilde{\eta}_1$ are
$O((\|\tilde{\eta}_1\|_2+\|H(\tilde{\eta})\|_2+\|\tilde{\eta}_3\|_2)^2(\|H(\tilde{\eta})\|_2+\|\tilde{\eta}_3\|_2)^2)$\linebreak
and hence $O(\mu\mu^{2+\alpha})=o(\mu^3),$ while
terms with three occurrences of $\tilde{\eta}_1$ are estimated by\linebreak
$O((\|\tilde{\eta}_1\|_{1,\infty} + \|\tilde{\eta}_{1xx}+k_0^2\tilde{\eta}_1\|_0)\|\tilde{\eta}_1\|_2^2(\|H(\tilde{\eta})\|_2+\|\tilde{\eta}_3\|_2)^2)
=O(\mu^{2+\alpha}\nn \tilde{\eta}_1\nn)=O(\mu^{5/2+\alpha})=o(\mu^3)$,
so that ${\mathcal K}_4(\tilde{\eta}) = -\frac{5}{4}\gamma\int_{-\infty}^\infty \tilde{\eta}_{1x}^2\tilde{\eta}_{1xx}^2+o(\mu^3)$.

Writing $\tilde{\eta}_1=\tilde{\eta}_1^++\tilde{\eta}_1^-$, where $\tilde{\eta}_1^+={\mathcal F}^{-1}[\chi_{[0,\infty)}{\mathcal F}[\tilde{\eta}_1]]$,
$\tilde{\eta}_1^-={\mathcal F}^{-1}[\chi_{(-\infty,0]}{\mathcal F}[\tilde{\eta}_1]]$, we find that
$$\|(\mathrm{i}k \mp \mathrm{i}k_0)\tilde{\eta}_1^\pm\|_s^2
=\|(|k|-k_0){\mathcal F}[\tilde{\eta}_1]\|_0^2
\leq \frac{1}{2}\int _{-\infty}^\infty (\mu^{2\alpha}+\mu^{-2\alpha}(|k|-k_0)^4) |{\mathcal F}[\tilde{\eta}_1]|^2\dk
\lesssim \mu^{2\alpha}\nn \tilde{\eta}_1 \nn^2
\lesssim \mu^{1+2\alpha}$$
so that $(\tilde{\eta}_1^\pm)_x = \pm \mathrm{i}k_0 + O(\mu^{1+2\alpha})$ in $H^s({\mathbb R})$ for each $s \geq 0$. Using this
estimate, one concludes that
\begin{align*}
\int_{-\infty}^\infty\tilde{\eta}_{1x}^2\tilde{\eta}_{1xx}^2 \dx &=\int_{-\infty}^\infty \left( (\tilde{\eta}_{1x}^+)^2(\tilde{\eta}_{1xx}^-)^2
+(\tilde{\eta}_{1x}^-)^2(\tilde{\eta}_{1xx}^+)^2+4\tilde{\eta}_{1x}\tilde{\eta}_{1x}^-\tilde{\eta}_{1xx}^+\tilde{\eta}_{1xx}^-\right)\dx \\
& =2k_0^6 \int_{-\infty}^\infty(\tilde{\eta}_1^+)^2(\tilde{\eta}_1^-)^2\dx+o(\mu)
\end{align*}
\begin{align*}
&=\frac{1}{3}k_0^2\int_{-\infty}^\infty\tilde{\eta}_1^4 \dx+ o(\mu).\qed
\end{align*}

The corresponding estimates for ${\mathcal L}_3(\tilde{\eta})$ and ${\mathcal L}_4(\tilde{\eta})$ are derived similarly
by Groves and Wahl\'{e}n \cite[\S4.3.2]{GrovesWahlen15}.

\begin{lemma}
Each near minimiser $\tilde{\eta}$ of ${\mathcal J}_\mu$ over $U \setminus \{0\}$ satisfies the estimates
$$- \nu_0^2{\mathcal L}_3(\tilde{\eta}) = A_3 \int_{-\infty}^\infty \tilde{\eta}_1^4 \dx + o(\mu^3), \qquad
{\mathcal L}_4(\tilde{\eta}_1)  = A_4^2 \int_{-\infty}^\infty \tilde{\eta}_1^4 \dx+ o(\mu^3).$$
\end{lemma}

\begin{corollary}
Suppose that $\tilde{\eta}$ is a near minimiser of ${\mathcal J}_\mu$ over $U \setminus \{0\}$. The estimates
\begin{align*}
{\mathcal M}_{a^2\mu}(a\tilde{\eta}) &= (a^3A_3+a^4A_4)\int_{-\infty}^\infty \tilde{\eta}_1^4 \dx+ a^3 o(\mu^3),\\
\langle {\mathcal M}_{a^2\mu}^\prime(a\tilde{\eta}), a\tilde{\eta} \rangle_0 + 4a^2\mu \tilde{{\mathcal M}}_{a^2\mu}(a\tilde{\eta})
&=(3a^3A_3+4a^4A_4)\int_{-\infty}^\infty \tilde{\eta}_1^4 \dx+ a^3 o(\mu^3)
\end{align*}
hold uniformly over $a \in [1,2]$, and $\int_{-\infty}^\infty \tilde{\eta}_1^4 \dx \gtrsim \mu^3$.
\end{corollary}

\begin{lemma} \label{SH step 1}
Suppose that $\tilde{\eta}$ is a near minimiser of ${\mathcal J}_\mu$ over $U \setminus \{0\}$ and $\frac{1}{2}A_3 + A_4<0$.
There exist $a_0 \in (1,2]$ and $q>2$ such that
$a \mapsto a^{-q}{\mathcal M}_{a^2\mu}(a\tilde{\eta})$, $a \in [1,a_0]$,
is decreasing and strictly negative.
\end{lemma}
{\bf Proof.} Observe that
\begin{eqnarray*}
\frac{\mathrm{d}}{\mathrm{d}a}\left(a^{-q}{\mathcal M}_{a^2\mu}(a\tilde{\eta})\right)
& = & a^{-(q+1)}\left(
-q{\mathcal M}_{a^2\mu}(a\tilde{\eta})+\langle {\mathcal M}_{a^2\mu}^\prime(a\tilde{\eta}),a\tilde{\eta}\rangle_0
+4a^2\mu\tilde{\mathcal M}_{a^2\mu}(a\tilde{\eta})
\right) \\
& = & a^{2-q}\left((3-q)A_3+a(4-q)A_4)\int_{\mathbb R} \tilde{\eta}_1^4 \dx + o(\mu^3)\right) \\
& \lesssim & -\mu^3 + o(\mu^3) \\
& <  & 0
\end{eqnarray*}
for  $a \in (1,a_0)$, $q \in (2,q_0)$; here $a_0>1$ and $q_0>2$ are chosen
so that $(3-q)A_3+a(4-q)A_4$, which is negative for $a=1$ and $q=2$,
is also negative for $a \in (1,a_0]$ and $q \in (2,q_0]$.\qed

\begin{corollary}
Suppose that $\frac{1}{2}A_3 + A_4<0$. The strict sub-homogeneity criterion $c_{a\mu} < ac_\mu$ holds for each $a>1$
(so that in particular $c_\mu$ is a strictly sub-additive function of $\mu$).
\end{corollary}
{\bf Proof.} It suffices to prove this inequality for $a \in (1,a_0^2]$. Let $\{\eta_n\}$ be a minimising sequence
for ${\mathcal J}_\mu$ over $U \setminus \{0\}$.
Replacing $a$ by $a^{1/2}$, we find from Lemma \ref{SH step 1} that
${\mathcal M}_{a\mu}(a^{1/2}\eta_n) \leq a^{1/2}q{\mathcal M}_\mu(\eta_n)$
and therefore that
$$
c_{a\mu} \leq {\mathcal J}_{a\mu}(\eta_n) \leq 
a\left({\mathcal K}_2(\eta_n)+\frac{\mu^2}{{\mathcal L}_2(\eta_n)}\right) + a^{1/2}q{\mathcal M}_{\mu}(\eta_n)
=a{\mathcal J}_\mu(\eta_n)+ (a^{1/2}q-a){\mathcal M}_\mu(\eta_n)$$
for $a \in (1,a_0^2]$. In the limit $n \rightarrow \infty$ this inequality yields
$c_{a\mu} < ac_\mu$ since $\limsup_{n \rightarrow \infty} {\mathcal M}_\mu(\eta_n)<0$.\qed

\begin{remark}
Theorem \ref{Convergence} is proved by Groves \& Wahl\'{e}n \cite[\S5.2.2]{GrovesWahlen15}; the
proof additionally confirms \emph{a posteriori} that the estimates $\nn \tilde{\eta}_1 \nn_\alpha^2 \lesssim \mu$,
$\|\tilde{\eta}_3\|_2^2 \lesssim\mu^{3+2\alpha}$ and $\|H(\tilde{\eta})\|_2^2 \lesssim\mu^{2+\alpha}$ also hold for $\alpha=1$.
\end{remark}
\noindent
{\bf Acknowledgement.} M. D. Groves would like to thank the Knut and Alice Wallenberg Foundation for funding a visiting professorship (reference KAW 2013.0318)
at Lund University during which this paper was prepared.

\end{document}